\documentclass[11pt]{amsart}

\usepackage{amsmath,amssymb}
\usepackage{mathtools}
\usepackage{amsthm}
\usepackage[abbrev]{amsrefs}
\usepackage{latexsym}
\usepackage{graphicx}
\usepackage{tikz}
\usepackage{bm}

\allowdisplaybreaks[1]

\title[]{$p$-adic Fourier theory for unramified extensions of $\mathbb{Q}_p$ and weighted Monna maps}
\date{}
\author{Yu Katagiri}

\newtheorem{thm}{Theorem}[section]
\newtheorem*{thm*}{Theorem}
\newtheorem{lem}[thm]{Lemma}
\newtheorem*{lem*}{Lemma}
\newtheorem{prop}[thm]{Proposition}
\newtheorem*{prop*}{Proposition}
\newtheorem{cor}[thm]{Corollary}
\newtheorem{cor*}{Corollary}

\theoremstyle{definition}
\newtheorem{dfn}[thm]{Definition}
\newtheorem*{dfn*}{Definition}
\newtheorem{ex}[thm]{Example}
\newtheorem*{ex*}{Example}
\newtheorem{rmk}[thm]{Remark}
\newtheorem*{rmk*}{Remark}

\newtheorem*{conj*}{Conjecture}

\newcommand{\F}{\mathbb{F}}
\newcommand{\Q}{\mathbb{Q}}
\newcommand{\Z}{\mathbb{Z}}
\newcommand{\R}{\mathbb{R}}
\newcommand{\C}{\mathbb{C}}

\newcommand{\G}{\mathcal{G}}

\newcommand{\LT}{\operatorname{LT}}

\makeatletter
\@addtoreset{equation}{section}

\makeatother

\makeatletter\@namedef{subjclassname@2020}{\textup{2020} Mathematics Subject Classification}\makeatother


\begin{document}
\maketitle

\begin{abstract}
We generalize a recent result of Ardakov and Berger on the $p$-adic valuations of Schneider--Teitelbaum polynomials evaluated at Lubin--Tate periods over the unramified quadratic extension of $\mathbb{Q}_p$. 
For finite unramified extensions of $\mathbb{Q}_p$ of degree $f \geq 2$, we introduce a weighted analogue of the Monna map and show that it determines the valuations of these polynomial values for Lubin--Tate formal groups whose formal logarithms are $p$-typical. 
Our proof uses the Artin--Hasse exponential and Witt vectors.
As an application, for a distinguished Lubin--Tate formal group, we determine the exact norms of the associated Schneider--Teitelbaum polynomials on Banach spaces of locally analytic functions and obtain an explicit rescaling of the Banach basis constructed by Bannai and Kobayashi so that each basis element has norm one.
\end{abstract}

\section{Introduction}

Let $p$ be a prime and let $\C_p$ denote the completion of an algebraic closure of $\Q_p$ with the absolute value $|\cdot |$ normalized by $|p|=p^{-1}$.
We regard every finite extension of $\Q_p$ as a subfield of $\C_p$.

Amice \cite{Am64} gave a description of locally analytic functions on $\Z_p$ in terms of the binomial polynomials $\binom{X}{n}$.
More precisely, for each positive integer $N$, the functions $\lfloor n/p^N \rfloor ! \binom{x}{n}$ form an orthonormal basis of the space of functions on $\Z_p$ that are analytic on each residue disk modulo $p^N$.

Schneider and Teitelbaum \cite{ST01} extended Amice's theory to arbitrary finite extensions of $\Q_p$ using Lubin--Tate formal groups.
Their theory, called $p$-adic Fourier theory, provides analogues of the classical binomial polynomials and the Amice transform.
We briefly describe the relevant objects below.

Let $K$ be a finite extension of $\Q_p$ with the residue field $\F_{q}$, where $q=p^f$, and the ramification index $e$.
We fix a uniformizer $\pi$ of $K$ and a Lubin--Tate formal group $\mathcal{G}$ attached to $\pi$.
Let $\log_{\mathcal{G}}(Z)$ denote its formal logarithm.
We choose a generator
\begin{align*}
G_\G(Z) \coloneqq \exp (\Omega_\G \log_{\mathcal{G}}(Z))-1
\end{align*}
of the rank-one $\mathcal{O}_K$-module $\operatorname{Hom}_{\mathcal{O}_{\C_p}}(\mathcal{G}, \widehat{\mathbb{G}_m})$.
The Lubin--Tate period $\Omega_\G$ is determined up to multiplication by an element in $\mathcal{O}_K^{\times}$.
It is also known that $\Omega_\G$ has the valuation
\begin{align*}
v_p(\Omega_\G)=\frac{1}{p-1}-\frac{1}{e(q-1)}.
\end{align*}
We define the {\it Schneider--Teitelbaum polynomials} $P_n^\G(X) \in K[X]$ by the generating series
\begin{align*}
\exp(X\log_{\mathcal{G}}(Z))=\sum_{n=0}^{\infty} P_n^\G(X)Z^n.
\end{align*}
The role of these polynomials in $p$-adic Fourier theory can be seen from the power series representation of the Fourier transform.
If $\mu$ is a locally $K$-analytic distribution on $\mathcal{O}_K$, then its Fourier transform can be expressed, in the coordinate associated with $\G$, as
\begin{align*}
\mathcal{F}(\mu)(Z)=\sum_{n=0}^{\infty} \mu(P_n^{\G}(x\Omega_{\G}))Z^n.
\end{align*}
Thus, the coefficients of the Fourier transform are obtained by evaluating distributions on the functions $x \mapsto P_n^{\G}(x\Omega_{\G})$.
(See \cite[Theorem 4.6]{BK16} or \cite[Theorems 2.3 and 3.6, and Lemma 4.6.9]{ST01}.)

When $K=\mathbb{Q}_p$, $\pi=p$ and $\G=\widehat{\mathbb{G}}_m$, we may take $\Omega_{\G}=1$, and the polynomials $P_n^\G(X)$ coincide with the classical binomial polynomials $\binom{X}{n}$. 
For a general Lubin--Tate formal group, the values $P_n^\G(x\Omega_\G)$ are $p$-adically integral whenever $x\in\mathcal{O}_K$. 
It is therefore natural to ask whether their $p$-adic valuations can be determined explicitly.

We now recall a recent result of Ardakov and Berger \cite{AB24b} on the valuations of Schneider--Teitelbaum polynomials.
Let $L$ be the unramified extension of $\Q_p$ of degree $f \geq 2$ and $q=p^f$, and take a uniformizer $\pi=p$.
In the present setting, we have
\begin{align*}
u \coloneqq v_p(\Omega_\G)=\frac{1}{p-1}-\frac{1}{q-1} = \frac{s}{q-1},
\end{align*}
where $s=\sum_{i=1}^{f-1} p^i$.
Let $\LT$ denote the Lubin--Tate formal group whose formal logarithm is
\begin{align}\label{log_LT}
\log_{\rm LT} (Z)=\sum_{i=0}^{\infty} \frac{Z^{q^i}}{p^i}.
\end{align}
%

\begin{thm}[{\cite[Theorem A]{AB24b}}]\label{thmAB}
Suppose that $f=2$.
For any nonnegative integer $n$, we have
\begin{align*}
v_p(P_n^{\LT}(\Omega_{\LT}))=w_2(n).
\end{align*}
Here, $w_2 : \Z_{\geq 0} \rightarrow \Q$ is given by
\begin{align}\label{monna}
w_2(n) \coloneqq u\sum_{i=0}^\infty n_ip^{-i}
\end{align}
if $n$ has the $p$-adic expansion $n=\sum_{i=0}^\infty n_ip^{i}$.
\end{thm}

The map sending $n$ to $p^{-1}u^{-1}w_2(n)$ is called the {\it Monna map}.

The aim of this paper is to generalize Theorem \ref{thmAB} to the case $f \geq 3$ and to other Lubin--Tate formal groups satisfying a certain condition.
As a generalization of \eqref{monna} to the case $f \geq 3$, we define $w_f : \Z_{\geq 0} \rightarrow \Q$ as the unique function satisfying the following conditions:
\begin{enumerate}
\item[(M1)] for any $0 \leq n \leq p^{f-1}-1$, we have $w_f(n)=nu-v_p(n!)$, 
\item[(M2)] for any $n \geq 0$ and $0 \leq i \leq p^{f-1}-1$, we have 
\begin{align*}
w_f(p^{f-1}n+i)=p^{-1}w_f(n)+w_f(i).
\end{align*}
\end{enumerate}
We note that it follows from \cite[Proposition 2.1]{AB24b} that \eqref{monna} satisfies the conditions (M1) and (M2) if $f=2$.
In this paper, we call $w_f$ the $f$-th {\it weighted Monna map}.
See Section 2 for properties of $w_f$, including a formula analogous to \eqref{monna} in terms of the $p$-adic expansion of $n$.

Recall that a formal power series of the form $\sum_{n=0}^\infty a_nZ^{p^n}$ is called {\it $p$-typical}.
Our main result is the following.

\begin{thm}\label{main1}
Let $f \geq 2$ and let $\G$ be a Lubin--Tate formal group over $\mathcal{O}_L$ attached to $p$ whose formal logarithm is $p$-typical.
%
%
Then, we have
\begin{align}\label{val=monna}
v_p(P_n^\G(\Omega_\G))=w_f(n)
\end{align}
for any $n \geq 0$.
\end{thm}

The argument used to prove Theorem \ref{thmAB} in \cite{AB24b} does not seem to extend directly to the case $f=3$, even when $\mathcal{G}=\LT$. 
Indeed, attempts to adapt the method of \cite[Section 4]{AB24b} lead to additional complications. 
We therefore take a different approach based on the Artin--Hasse exponential and the theory of Witt vectors, which also enables us to treat Lubin--Tate formal groups with $p$-typical formal logarithms.

We also note that, if the formal logarithm of $\G$ is not $p$-typical, then \eqref{val=monna} does not hold in general.
(See Remark \ref{bad_ex}.)
However, if $n$ is a $p$-power, then \eqref{val=monna} holds for any $\G$.
(See Proposition \ref{thm_for_powers}.)

Theorem \ref{main1} also leads to an application in $p$-adic functional analysis.
In their study of integral structures in $p$-adic Fourier theory, Bannai and Kobayashi \cite{BK16} constructed Banach bases of locally analytic function spaces using suitably normalized Schneider--Teitelbaum polynomials.
They also obtained uniform estimates for the norms of these basis elements.
In \cite{Ka26}, these estimates were refined to exact norm formulas for the distinguished Lubin--Tate formal group $\LT$ over the unramified quadratic extension of $\Q_p$.
We show that the valuation formula in Theorem \ref{main1} allows us to extend this norm computation to arbitrary unramified extensions of finite degree.

For a positive integer $N$, let $\mathrm{LA}_N(\mathcal{O}_L,\mathbb{C}_p)$ denote the space of functions on $\mathcal{O}_L$ that are analytic on each residue disk $a+p^N\mathcal{O}_L$. 
For such a function $f$, write
\begin{align*}
f(a+p^Nx)=\sum_{j=0}^{\infty}c_{a,j}x^j,
\end{align*}
where $c_{a, j} \to 0$ as $j \to \infty$, for each $a\in\mathcal{O}_L$. 
We equip this space with the norm
\begin{align*}
||f||_N \coloneqq \sup_{a\in\mathcal{O}_L}\sup_{j\geq0}\{|c_{a,j}|\}.
\end{align*}
Our application concerns the exact values of this norm for the functions $x\mapsto P_n^{\LT}(x\Omega_{\LT})$.

\begin{thm}\label{cor1}
Let $N$ be a positive integer.
For $n \geq 0$, write $n=lq^{N+1}+r$ with $l \geq 0$ and $0 \leq r \leq q^{N+1}-1$.
Then we have
\begin{align*}
|| P_n^{\LT} (x\Omega_{\LT})||_N =p^{-\eta_N(n)},
\end{align*}
where
\begin{align*}
\eta_N(n) \coloneqq w_f(r)+l(u-1)-v_p(l!).
\end{align*}
\end{thm}

The explicit norm formula also provides a refinement of the Bannai--Kobayashi normalization.
More precisely, Theorems \ref{main1} and \ref{cor1} allow us to rescale their Banach basis so that all its elements have norm one.

\begin{cor}\label{cor2}
For  $n=lq^{N+1}+r$ with $l \geq 0$ and $0 \leq r \leq q^{N+1}-1$, we define
\begin{align*}
\widetilde{e_{N, n}}(x) \coloneqq 
P_r^{\LT}(\Omega_{\LT})^{-1} \underline{\gamma} \left( \left\lfloor \frac{n}{q^N} \right\rfloor \right) P_n^{\LT}(x\Omega_{\LT}).
\end{align*}
Then, we have $|| \widetilde{e_{N, n}}(x) ||_N=1$. 
(See Section 5 for the definition of $\underline{\gamma}$.)
\end{cor}

The paper is organized as follows. 
In Section 2, we study the weighted Monna map $w_f$, obtaining an explicit formula in terms of the $p$-adic expansion and several inequalities that will be used in the proof of Theorem \ref{main1}.
In Section 3, we recall some preliminary results on Lubin--Tate formal groups, Witt vectors, and the Artin--Hasse exponential.
In Section 4, we prove Theorem \ref{main1} using the Artin–Hasse method and establish a valuation formula at powers of $p$ for arbitrary Lubin--Tate formal groups.
Finally, in Section 5, we prove Theorem \ref{cor1} and Corollary \ref{cor2} by adapting the arguments of \cite{Ka26}.

\vspace{10pt}
\noindent
\textsc{Notation:}~ For a nonnegative integer $a$, we write the $p$-adic expansion of $a$ by
\begin{align*}
a=\sum_{i=0}^{\infty} a_ip^i
\end{align*}
if there is no confusion.
We use the convention $v_p(0)=+\infty$.

\vspace{10pt}
\noindent
\textsc{Acknowledgments:}~ This work started from Shinichi Kobayashi's comments on the author's work \cite{Ka26}. 
The author is grateful to him for many helpful comments and discussion.
The author is also grateful to Laurent Berger for his helpful comments and stimulating questions.
The author would like to thank Ryota Shii, Masaya Miyazaki and Shoma Sueyoshi for helpful discussion.

\section{The map $w_f$}

In this section, we establish some properties of the weighted Monna map $w_f$ that will be used in the proofs of our main results.

\subsection{Basic properties}

We first derive an explicit formula for $w_f$ in terms of the $p$-adic expansion of its argument and establish some basic properties.

\begin{lem}\label{w_at_p}
For any $j \geq 0$ and $0 \leq r \leq f-2$, we have
\begin{align*}
w_f(p^{(f-1)j+r})=\frac{s+1-p^r}{p^j(q-1)}<1.
\end{align*}
\end{lem}

\begin{proof}
We check the assertion by induction on $j$.
For $j=0$, since $p^r \leq p^{f-1}-1$, the condition (M1) implies that
\begin{align*}
w_f(p^r)&=p^ru-v_p(p^r!) \\
&=p^r \left( \frac{1}{p-1}-\frac{1}{q-1} \right) -\frac{p^r-1}{p-1} \\
&=\frac{1}{p-1}-\frac{p^r}{q-1}=\frac{s+1-p^r}{q-1}.
\end{align*}
We also have
\begin{align*}
\frac{s+1-p^r}{q-1}<\frac{s+1}{q-1} \cdot (p-1)=1.
\end{align*}

Suppose that $j>0$ and that the lemma holds for $j-1$.
Then it follows from the condition (M2) and induction hypothesis that
\begin{align*}
w_f(p^{(f-1)j+r})=p^{-1} w_f(p^{(f-1)(j-1)+r})=\frac{s+1-p^r}{p^j(q-1)}<1.
\end{align*}
\end{proof}

\begin{prop}\label{explicit_w_f}
If a nonnegative integer $n$ has the $p$-adic expansion $n=\sum_{i=0}^{\infty} n_ip^i$, then we have
\begin{align}\label{rewrite_w_f}
w_f(n)=\sum_{i=0}^{\infty} n_iw_f(p^i) = \sum_{j=0}^{\infty} \sum_{r=0}^{f-2} \frac{s+1-p^r}{p^j(q-1)} n_{(f-1)j+r}.
\end{align}
\end{prop}

\begin{proof}
The second equality follows from Lemma \ref{w_at_p}, and hence we prove the first equality by induction on $n$.
If $0 \leq n \leq p^{f-1}-1$, the condition (M1) and Lemma \ref{w_at_p} give
\begin{align*}
w_f(n)&=nu-v_p(n!)=nu-\frac{1}{p-1}n+\frac{1}{p-1}\sum_{i=0}^{f-2} n_i \\
&=-\frac{1}{q-1}\sum_{i=0}^{f-2} n_i p^i+\frac{s+1}{q-1}\sum_{i=0}^{f-2} n_i \\
&=\sum_{i=0}^{f-2} \frac{s+1-p^i}{q-1}n_i=\sum_{i=0}^{f-2} n_iw_f(p^i).
\end{align*}

Suppose that $n \geq p^{f-1}$ and that the proposition holds for all nonnegative integers less than $n$.
Since
\begin{align*}
n=p^{f-1}\sum_{i=0}^{\infty} n_{i+f-1}p^i+\sum_{i=0}^{f-2} n_i p^i
\end{align*}
and the first summation on the right-hand side is less than $n$, the condition (M2) and the induction hypothesis imply that 
\begin{align*}
w_f(n)&=p^{-1}w_f \left( \sum_{i=0}^{\infty} n_{i+f-1}p^i \right)+w_f \left(\sum_{i=0}^{f-2} n_i p^i \right) \\
&=p^{-1} \sum_{i=0}^\infty n_{i+f-1}w_f(p^i)+\sum_{i=0}^{f-2} n_iw_f(p^i) \\
&=\sum_{i=0}^\infty n_{i+f-1}w_f(p^{i+f-1})+\sum_{i=0}^{f-2} n_iw_f(p^i)=\sum_{i=0}^{\infty} n_iw_f(p^i).
\end{align*}
\end{proof}

The following proposition gives some basic properties of $w_f$, which are not needed for the proof of Theorem \ref{main1}. 
Some of them will be used in the proof of Corollary \ref{comm_main1} and in Section 5.

\begin{prop}[cf. {\cite[Proposition 2.1]{AB24b}}]\label{prop2.1}
The function $w_f$ has the following properties.
\begin{enumerate}
\item For $m \geq 1$, $n \geq 0$ and $0 \leq i \leq p^m-1$, we have
\begin{align*}
w_f(p^mn+i)=w_f(p^mn)+w_f(i).
\end{align*}
\item For a nonnegative integer $n$ whose $p$-adic expansion $n=\sum_{i=0}^{\infty} n_ip^i$, we define
\begin{align*}
D_f(n) \coloneqq \sum_{j=0}^{\infty} \left(\sum_{r=0}^{f-2} n_{(f-1)j+r} \right){p^{-j-1}}.
\end{align*}
For $n \geq 0$, we have $w_f(n)-w_f(qn)=D_f(n)$.
\item For $n \geq 0$, we have
\begin{align*}
w_f(n)=\sum_{i=0}^{\infty} D_f(q^in).
\end{align*}
\item For all $a, b \geq 0$, we have $w_f(a+b) \leq w_f(a)+w_f(b)$.
\item We have
\begin{align*}
\sup_{n \geq 0} \{ w_f(n) \}=f-1+\frac{f-2}{p-1}+\frac{1}{q-1}.
\end{align*}
\end{enumerate}
\end{prop}

\begin{proof}
Properties (1) and (2) follow from a direct calculation.
Property (3) follows from (2) since $w_f(q^in) \to 0$ as $i \to \infty$.
We prove (4).
From (3), it suffices to show that
\begin{align*}
D_f(a+b) \leq D_f(a)+D_f(b).
\end{align*}
We put $c \coloneqq a+b$ and define the sequence $(d_i)_{i \geq -1}$ of $\{ 0, 1\}$ by
\begin{align*}
d_{-1}=0, \qquad c_{i}=a_{i}+b_{i}+d_{i-1}-pd_{i}.
\end{align*}
If we write $i=(f-1)j+r$ with $j \geq 0$ and $0 \leq r \leq f-2$, and put $\lambda_i \coloneqq p^{-j-1}$, then we have
\begin{align*}
&D_f(a)+D_f(b)-D_f(c) \\
=&\sum_{i=0}^{\infty} (pd_i-d_{i-1})\lambda_i=\sum_{i=0}^{\infty} (p\lambda_i-\lambda_{i+1})d_i \geq 0.
\end{align*}

Finally, since $0 \leq n_i \leq p-1$ for all $i \geq 0$, it follows from \eqref{rewrite_w_f} that
\begin{align*}
w_f(n) &\leq (p-1) \sum_{j=0}^{\infty} \sum_{r=0}^{f-2} \frac{s+1-p^r}{p^j(q-1)} \\
&=\frac{p}{q-1} \sum_{r=0}^{f-2} (s+1-p^r) =f-1+\frac{f-2}{p-1}+\frac{1}{q-1}.
\end{align*}
On the other hand, by \eqref{rewrite_w_f}, we have
\begin{align*}
\lim_{m \to \infty} w_f(p^m-1)=(p-1) \sum_{i=0}^{\infty} w_f(p^i)=f-1+\frac{f-2}{p-1}+\frac{1}{q-1}.
\end{align*}
This proves (5).
\end{proof}

\begin{rmk}
Proposition \ref{prop2.1} gives another characterization of $w_f$, that is, the function $w_f$ is the unique map $w_f : \Z_{\geq 0} \rightarrow \R$ satisfying 
\begin{enumerate}
\item for any $n \geq 0$, we have $w_f(p^{f-1}n)=p^{-1}w_f(n)$,
\item for any $n \geq 0$, we have $w_f(n)-w_f(qn)=D_f(n)$.
\end{enumerate}
\end{rmk}

\subsection{Inequalities for $w_f$}

We next study the values of $w_f$ at powers of $p$. 
The resulting inequalities will play a key role in the proof of Theorem \ref{main1}.

\begin{lem}\label{diff_p}
For any $m \geq 0$, we have
\begin{align*}
w_f(p^m)>w_f(p^{m+1}).
\end{align*}
Furthermore, for any $1 \leq n \leq p^m-1$, we have
\begin{align*}
w_f(n) > w_f(p^m).
\end{align*}
\end{lem}

\begin{proof}
We write $m=(f-1)j+r$ with $j \geq 0$ and $0 \leq r \leq f-2$. 
If $r=f-2$, we have
\begin{align*}
&w_f(p^m)-w_f(p^{m+1}) \\
=& \frac{s+1-p^{f-2}}{p^{j}(q-1)}-\frac{s}{p^{j+1}(q-1)} \\
=&\frac{1}{p^{j+1}(q-1)} \left\{ p(s+1-p^{f-2})-s \right\} =\frac{p^{f-2}(p-1)}{p^{j}(q-1)}>0.
\end{align*}
For $0 \leq r \leq f-3$, we similarly have
\begin{align*}
w_f(p^m)-w_f(p^{m+1}) =\frac{p^{r}(p-1)}{p^{j}(q-1)}>0.
\end{align*}
This proves the first assertion.

The second assertion follows from \eqref{rewrite_w_f} since
\begin{align*}
w_f(n)=\sum_{j=0}^{m-1} n_jw_f(p^j) \geq w_f(p^{m-1}) >w_f(p^m).
\end{align*}
\end{proof}

\begin{lem}\label{carrylem}
Let $n$ be a positive integer whose  $p$-adic expansion $n=\sum_{i=0}^{\infty} n_ip^i$ and let $\nu=(\nu_i)_{i \geq 0}$ be a sequence of nonnegative integers with finite support satisfying $\sum_{i=0}^{\infty} \nu_ip^i=n$.
\begin{enumerate}
\item We define the sequence $(c_i)_{i \geq -1}$ of integers by $c_{-1}=0$ and 
\begin{align*}
\nu_i+c_{i-1}=n_i+pc_i
\end{align*}
for $i \geq 0$. 
Then, we have
\begin{align}\label{expr_c}
p^{i+1}c_i=\sum_{l=0}^i (\nu_l-n_l)p^l
\end{align}
and $c_i \geq 0$ for any $i \geq 0$.
\item If $\nu \neq (n_i)_{i \geq 0}$, then we have
\begin{align*}
\sum_{i=0}^{\infty} \nu_i w_f(p^i)>w_f(n).
\end{align*}
\end{enumerate}
\end{lem}

\begin{proof}
\begin{enumerate}
\item Since $p^{l+1}c_l-p^lc_{l-1}=(\nu_l-n_l)p^l$ by the definition of $(c_i)_{i \geq -1}$, we obtain \eqref{expr_c}.
On the other hand, by the assumption, we have
\begin{align}\label{nu_and_residue}
\sum_{l=0}^i \nu_lp^l \equiv \sum_{l=0}^i n_lp^l \bmod p^{i+1}.
\end{align}
Since the right-hand side of \eqref{nu_and_residue} is the remainder of $n$ modulo $p^{i+1}$ and the left-hand side is nonnegative, we find that $p^{i+1}c_i \geq 0$ by \eqref{expr_c}.
\item By the definition of $(c_i)_{i \geq -1}$ and \eqref{rewrite_w_f}, we find that 
\begin{align*}
\sum_{i=0}^{\infty} \nu_i w_f(p^i)&=\sum_{i=0}^{\infty} (n_i+pc_i-c_{i-1})w_f(p^i) \\
&=w_f(n)+\sum_{i=0}^{\infty} (pw_f(p^i)-w_f(p^{i+1}))c_i.
\end{align*}
If $\nu \neq (n_i)_{i \geq 0}$, there exists $i \geq 0$ such that $c_i>0$.
Hence Lemma \ref{diff_p} implies the claim.
\end{enumerate}
\end{proof}

\section{Preliminaries}

In this section, we collect some preliminary results on Lubin–Tate formal groups, Witt vectors, and the Artin–Hasse exponential.

\subsection{Basic properties of Lubin--Tate formal groups}

We first recall two basic properties of Lubin--Tate formal groups that will be used in the proof of Theorem \ref{main1}.
Let $\mathcal{G}$ be a Lubin--Tate formal group attached to $p$.
In this subsection, we do not assume that the formal logarithm of $\G$ is $p$-typical.
If there is no confusion, we suppress the superscript or subscript $\G$ in $P_n^\G$, $\Omega_\G$ and $G_\G$.
Put
\begin{align*}
G(Z) \coloneqq \exp (\Omega \log_{\mathcal{G}}(Z))-1=\sum_{n=1}^\infty P_n(\Omega)Z^n.
\end{align*}

\begin{lem}\label{ABsection3}
We have
\begin{align*}
G(Z^q) \equiv G(Z)^p \bmod p\mathcal{O}_{\C_p}.
\end{align*}
\end{lem}

\begin{proof}
The case for $f=2$ and $\G=\LT$ follows from \cite[Proposition 3.1 and Corollary 3.4]{AB24b}.
One can show the lemma for a general case by the same argument.
\end{proof}

\begin{lem}\label{int_coeff_log}
For any $1 \leq n \leq q-1$, the coefficient of $Z^n$ in $\log_\G(Z)$ is integral.
\end{lem}

\begin{proof}
There exists $\psi(Z) \in Z+Z^2\mathcal{O}_L[[Z]]$ such that $\log_\G=\log_{\LT} \circ \psi$.
Since $\log_{\LT}$ is given by \eqref{log_LT}, we find that $\log_\G(Z) \equiv \psi(Z) \bmod Z^q$.
\end{proof}

\subsection{Witt vectors and the Artin--Hasse exponential}
The Artin--Hasse exponential and the theory of $p$-typical Witt vectors provide the main tools for the proof of Theorem \ref{main1}.
We recall the necessary definitions and properties, following \cite[Chapter I\hspace{-1.2pt}I \S 6]{Se79} and \cite[Chapter I \S8-9]{FV02}.

Let $A$ be a unital commutative ring and let $W(A)$ denote the ring of $p$-typical Witt vectors over $A$.
As a set, $W(A)$ consists of sequences $(x_0, x_1, \cdots)$ with $x_i \in A$, although its ring operations are not defined componentwise.
We define {\it Verschiebung} $V$ on $W(A)$ by
\begin{align*}
V\bm{x} \coloneqq (0, x_0, x_1, \cdots, x_{i-1}, \cdots) \in W(A)
\end{align*}
for $\bm{x}=(x_i)_{i \geq 0} \in W(A)$.
It is known that $V$ is an injective additive map.

In addition, if $A$ is a ring of characteristic $p$, we also define a ring endmorphism {\it Frobenius} $F$ on $W(A)$ by
\begin{align*}
F\bm{x} \coloneqq (x_0^p, x_1^p, \cdots, x_{i}^p, \cdots) \in W(A)
\end{align*}
for $\bm{x}=(x_i)_{i \geq 0} \in W(A)$.
Furthermore, we have an important property
\begin{align}\label{FV}
FV=VF=p,
\end{align}
where $p$ on the right-hand side of \eqref{FV} means multiplication by $p$ in $W(A)$.

We define the {\it Artin--Hasse exponential} $E_p(X) \in 1+X\Q[[X]]$ by
\begin{align*}
E_p(X) \coloneqq \exp \left( \sum_{i=0}^{\infty} \frac{X^{p^i}}{p^i} \right).
\end{align*}
Let $\alpha_n$ denote the coefficient of $X^n$ in $E_p(X)$.
The following is a fundamental property of $E_p$.

\begin{prop}\label{intE_p}
We have $E_p(X) \in 1+X\Z_p[[X]]$.
Moreover, if $0 \leq n \leq p-1$, we have
\begin{align*}
\alpha_n=\frac{1}{n!}.
\end{align*}
\end{prop}

\begin{proof}
The first assertion is well-known.
(See \cite[pp.92-95]{Ko84} for example.)
The second assertion follows from the fact that
\begin{align*}
E_p(X) \equiv \exp (X) \bmod X^p.
\end{align*}
\end{proof}

Since $E_p(X)$ has rational coefficients, Proposition \ref{intE_p} implies that $E_p(X) \in 1+X\Z_{(p)}[[X]]$.
Hence, we can define the following map for any $\Z_{(p)}$-algebra.

\begin{dfn}[{cf. \cite[Chapter I \S9]{FV02}}]
Let $A$ be a $\Z_{(p)}$-algebra.
We define the {\it Artin--Hasse map} $\mathcal{E}_A : W(A) \rightarrow 1+ZA[[Z]]$ by
\begin{align*}
\mathcal{E}_A(\bm{x}; Z) \coloneqq \prod_{i=0}^\infty E_p(x_iZ^{p^i})
\end{align*}
for $\bm{x}=(x_i)_{i \geq 0} \in W(A)$.
\end{dfn}

The map $\mathcal{E}_A$ has the following fundamental properties.

\begin{prop}[{\cite[Chapter I \S9]{FV02}}]\label{propertyE_A}
\begin{enumerate}
\item For any $\bm{x}, \bm{y} \in W(A)$, we have $\mathcal{E}_A(\bm{x}+\bm{y}; Z)=\mathcal{E}_A(\bm{x}; Z)\mathcal{E}_A(\bm{y}; Z)$.
\item For any $\bm{x} \in W(A)$, we have $\mathcal{E}_A(V\bm{x}; Z)=\mathcal{E}_A(\bm{x}; Z^p)$.
\item The map $\mathcal{E}_A$ is injective.
\end{enumerate}
\end{prop}

We will use the following classical lemma concerning the Artin--Hasse exponential.
(See, for example, \cite[Chapter 3]{De72} or \cite[Subsection 3.2]{Ri16}.)

\begin{lem}\label{Wittlem}
Let $f(Z) \in 1+Z\C_p[[Z]]$ satisfy
\begin{align*}
\log f(Z)=\sum_{n=0}^{\infty} a_n Z^{p^n}.
\end{align*}
We set
\begin{align*}
b_0=a_0, \qquad b_n \coloneqq a_n-\sum_{i=0}^{n-1} \frac{b_i^{p^{n-i}}}{p^{n-i}}
\end{align*}
and $\bm{b}=(b_i)_{i \geq 0} \in W(\C_p)$.
\begin{enumerate}
\item We have $f(Z)=\mathcal{E}_{\C_p}(\bm{b}; Z)$.
%
%
\item If $f(Z) \in 1+Z\mathcal{O}_{\C_p}[[Z]]$, then we have $b_n \in \mathcal{O}_{\C_p}$ for all $n \geq 0$.
\end{enumerate}
\end{lem}

\begin{proof}
\begin{enumerate}
\item 
By the recursive definition of $b_n$, we have
\begin{align*}
\log \prod_{n=0}^{\infty} E_p(b_nZ^{p^n}) &=\sum_{n=0}^{\infty} \sum_{i=0}^{\infty} \frac{b_n^{p^i}}{p^i}Z^{p^{n+i}}=\sum_{n=0}^{\infty} \left(\sum_{i=0}^n \frac{b_i^{p^{n-i}}}{p^{n-i}} \right)Z^{p^n} \\
&=\sum_{n=0}^{\infty} a_nZ^{p^n}=\log f(Z).
\end{align*}
\item We show the assertion by induction on $n$.
The case for $n=0$ follows from $b_0=a_0$.

Suppose that $n>0$ and that $b_0, \cdots, b_{n-1} \in \mathcal{O}_{\C_p}$.
Then, since $\prod_{i=0}^{n-1} E_p(b_iZ^{p^i}) \in 1+Z\mathcal{O}_{\C_p}[[Z]]$ by induction hypothesis and Proposition \ref{intE_p}, we have
\begin{align*}
f(Z)\left(\prod_{i=0}^{n-1} E_p(b_iZ^{p^i})\right)^{-1} =\prod_{i=n}^{\infty} E_p(b_iZ^{p^i})=1+b_nZ^{p^n}+O(Z^{p^n+1})
\end{align*}
by the assertion (1).
Hence $b_n \in \mathcal{O}_{\C_p}$.
\end{enumerate}
\end{proof}

\section{Proof of Theorem \ref{main1}}

\subsection{The Artin--Hasse method}

In this subsection, we suppose that the formal logarithm of $\G$ is $p$-typical and is given by
\begin{align*}
\log_\G (Z)=\sum_{n=0}^{\infty} a_nZ^{p^n}.
\end{align*}
Hence, if we put
\begin{align*}
f(Z) \coloneqq 1+G(Z)= \exp( \Omega \log_{\G} (Z))=\sum_{n=0}^{\infty} P_n(\Omega)Z^n,
\end{align*}
then we have
\begin{align*}
\log f(Z)=\Omega \log_{\G} (Z) = \sum_{n=0}^{\infty} \Omega a_nZ^{p^{n}}.
\end{align*}
We now set
\begin{align*}
b_0 \coloneqq \Omega a_0=\Omega, \qquad b_n \coloneqq \Omega a_n-\sum_{i=0}^{n-1} \frac{b_i^{p^{n-i}}}{p^{n-i}}.
\end{align*}
Applying Lemma \ref{Wittlem}(1), we obtain
\begin{align*}
\sum_{n=0}^{\infty} P_n(\Omega)Z^n&=\prod_{i=0}^{\infty} E_p(b_iZ^{p^i})=\prod_{i=0}^{\infty} \sum_{\nu_i=0}^{\infty} \alpha_{\nu_i}(b_iZ^{p^i})^{\nu_i},
\end{align*}
and hence
\begin{align}\label{exp_u_n}
P_n(\Omega)=\sum_{\nu=(\nu_i)_{i}} \prod_{i=0}^{\infty} \alpha_{\nu_i}b_i^{\nu_i}.
\end{align}
Here, $\nu$ in the summation in \eqref{exp_u_n} runs through all sequences $\nu=(\nu_i)_{i \geq 0}$ of nonnegative integers with finite support satisfying $\sum_{i=0}^{\infty} \nu_ip^i=n$.
We also have $\bm{b}=(b_i)_{i \geq 0} \in W(\mathcal{O}_{\C_p})$ by Lemma \ref{Wittlem}(2).

To apply the expression \eqref{exp_u_n} to the proof of Theorem \ref{main1}, we determine the $p$-adic valuation of $b_n$ for each $n \geq 0$.
We put $A \coloneqq \mathcal{O}_{\C_p}/{p\mathcal{O}_{\C_p}}$ and let $\overline{\bm{b}}$ denote the image of $\bm{b}$ under the natural map $W(\mathcal{O}_{\C_p}) \rightarrow W(A)$.

\begin{lem}\label{conseq_Witt}
We have 
\begin{align}\label{rel_VandF}
V^{f-1}\overline{\bm{b}}=F\overline{\bm{b}} \in W(A).
\end{align}
In particular, we have
\begin{align*}
b_{n+f-1}^p \equiv b_n \bmod p\mathcal{O}_{\C_p}
\end{align*}
for any $n \geq0$.
\end{lem}

\begin{proof}
Since $\bm{b}=(b_i)_{i \geq 0} \in W(\mathcal{O}_{\C_p})$ and $E_p(X) \in 1+X\Z_p[[X]]$ by Proposition \ref{intE_p}, we see that 
\begin{align*}
1+G(Z) =\prod_{i=0}^{\infty} E_p(b_iZ^{p^i}) \equiv \mathcal{E}_A(\overline{\bm{b}}; Z) \bmod p\mathcal{O}_{\C_p}.
\end{align*}
Then, Proposition \ref{propertyE_A} and Lemma \ref{ABsection3} imply that
\begin{align*}
\mathcal{E}_A(V^f\overline{\bm{b}}; Z) &=\mathcal{E}_A(\overline{\bm{b}}; Z^q) \equiv 1+G(Z^q) \\
&\equiv 1+G(Z)^p \equiv \mathcal{E}_A(p\overline{\bm{b}}; Z) \bmod p\mathcal{O}_{\C_p}.
\end{align*}
Since $\mathcal{E}_A : W(A) \rightarrow 1+ZA[[Z]]$ is injective by Proposition \ref{propertyE_A}(3), we obtain $V^f\overline{\bm{b}}=p\overline{\bm{b}}=VF\overline{\bm{b}} \in W(A)$ by \eqref{FV}.
Since $V$ is also injective, \eqref{rel_VandF} follows.

Comparing the $(n+f-1)$-th components on both sides of \eqref{rel_VandF}, we obtain the second assertion.
\end{proof}

\begin{lem}\label{n=p^m}
For any $n \geq 0$, we have
\begin{align*}
v_p(b_n)=w_f(p^n).
\end{align*}
\end{lem}

\begin{proof}
We prove the lemma by induction on $n$.
For $n=0$, it follows from the fact that $b_0=P_1(\Omega)=\Omega$.
Suppose that $n>0$ and that the claim is true for all nonnegative integers less than $n$.

We first assume that $n \leq f-2$.
For $0 \leq i \leq n-1$, the induction hypothesis and Lemma \ref{w_at_p} imply that
\begin{align*}
v_p(b_i^{p^{n-i}}/p^{n-i})&=p^{n-i}w_f(p^i)-(n-i) \\
&=p^{n-i} \left( \frac{1}{p-1}-\frac{p^i}{q-1} \right)-(n-i) \\
&=\left( \frac{1}{p-1}-\frac{p^n}{q-1} \right)+\frac{p^{n-i}-1}{p-1}-(n-i) \\
&=w_f(p^n)+\sum_{k=1}^{n-i}(p^{k-1}-1) \geq w_f(p^n).
\end{align*}
This also shows that $v_p(b_i^{p^{n-i}}/p^{n-i})=w_f(p^n)$ if and only if $i=n-1$.
Since $w_f(p^n)<w_f(1)=u$ by Lemma \ref{diff_p}, and $a_n \in \mathcal{O}_L$ by Lemma \ref{int_coeff_log}, the definition of $b_n$ gives
\begin{align*}
v_p(b_n)=\min_{0 \leq i \leq n-1} \left\{ v_p(\Omega a_n), v_p(b_i^{p^{n-i}}/p^{n-i}) \right\}=w_f(p^n).
\end{align*}

We finally consider the case $n \geq f-1$.
Since $v_p(b_{n-f+1})=w_f(p^{n-f+1})<1$ by the induction hypothesis and Lemma \ref{w_at_p}, we find that
\begin{align*}
pv_p(b_n)&=v_p(b_n^p-b_{n-f+1}+b_{n-f+1}) \\
&=\min \{ v_p(b_n^p-b_{n-f+1}), v_p(b_{n-f+1}) \}=v_p(b_{n-f+1})=w_f(p^{n-f+1}).
\end{align*}
Hence, we obtain
\begin{align*}
v_p(b_n)=p^{-1}w_f(p^{n-f+1})=w_f(p^n)
\end{align*}
by the condition (M2) in the definition of $w_f$.
\end{proof}

We prove Theorem \ref{main1}.

\begin{proof}[Proof of Theorem \ref{main1}]
Let $n=\sum_{i=0}^\infty n_ip^i$ denote the $p$-adic expansion of a nonnegative integer $n$.
If we let $\nu'=(n_i)_{i \geq 0}$, \eqref{exp_u_n} again gives
\begin{align*}
P_{n}(\Omega)=\prod_{i=0}^{\infty} \frac{b_i^{n_i}}{n_i!}+\sum_{\nu\neq \nu'} \prod_{i=0}^{\infty} \alpha_{\nu_i}b_i^{\nu_i}.
\end{align*}
Here, $\nu$ in the summation runs through all sequences $\nu=(\nu_i)_{i \geq 0}$ of nonnegative integers with finite support satisfying $\sum_{i=0}^{\infty} \nu_ip^i=n$ and $\nu\neq \nu'$, and we use the fact that $0 \leq n_i \leq p-1$ and Proposition \ref{intE_p}.
For $\nu \neq \nu'$, it follows from Proposition \ref{intE_p} and Lemma \ref{n=p^m} that
\begin{align*}
v_p \left( \prod_{i=0}^{\infty} \alpha_{\nu_i}b_i^{\nu_i} \right) \geq \sum_{i=0}^{\infty} \nu_i v_p(b_i)=\sum_{i=0}^{\infty} \nu_i w_f(p^i).
\end{align*}
Since $\nu'$ corresponds to the $p$-adic expansion of $n$, Lemma \ref{carrylem} implies that
\begin{align*}
\sum_{i=0}^{\infty} \nu_i w_f(p^i)>w_f(n).
\end{align*}
By \eqref{rewrite_w_f}, we also find that 
\begin{align*}
v_p \left( \prod_{i=0}^{\infty} \frac{b_i^{n_i}}{n_i!} \right) = \sum_{i=0}^{\infty} n_i v_p(b_i)=\sum_{i=0}^{\infty} n_i w_f(p^i)=w_f(n).
\end{align*}
Thus, we obtain $v_p(P_n(\Omega))=w_f(n)$.
\end{proof}

\begin{ex}
The following is due to Hazewinkel \cite[Section 5]{Ha79}.
(See also \cite[Chapter I]{Ha78}.)
Let $g(Z) \in Z+Z^2\mathcal{O}_L[[Z]]$ be a $p$-typical formal power series and put
\begin{align*}
f(Z)=\sum_{n=0}^\infty \frac{g(Z^{q^n})}{p^n}.
\end{align*}
Then, $\G(X, Y)=f^{-1}(f(X)+f(Y))$ is a Lubin--Tate formal group attached to $p$ and its formal logarithm is given by $f$, which is $p$-typical.
\end{ex}

\begin{rmk}\label{bad_ex}
Let $p$ be an odd prime and $f \geq 2$.
The valuation formula \eqref{val=monna} does not hold for arbitrary Lubin--Tate formal groups.
Indeed, if $[p]_\G(Z)=pZ+pZ^2+Z^q$, then we find that
\begin{align*}
\log_{\G}(Z)=Z+\frac{1}{1-p}Z^2+O(Z^3)
\end{align*}
and that
\begin{align*}
P_2^\G(X)=\frac{1}{2}X^2+\frac{1}{1-p}X.
\end{align*}
Then we have
\begin{align*}
v_p(P_2^\G(\Omega_\G))=\min \{ 2u, u\}=u<2u=w_f(2).
\end{align*}
\end{rmk}

We will use the following corollary in Section 5.

\begin{cor}[cf. {\cite[Lemma 2.2]{Ka26}}]\label{comm_main1}
Let $\G$ be a Lubin--Tate formal group attached to $p$ whose formal logarithm is $p$-typical.
For any $n \geq 0$ and $a \in \mathcal{O}_L$, we have
\begin{align*}
|P_n(a\Omega)| \leq p^{-w_f(n)}.
\end{align*}
Moreover, equality holds if $a \in \mathcal{O}_L^{\times}$.
\end{cor}

\begin{proof}
This follows from the same argument as in \cite[Lemma 2.2]{Ka26}.
\end{proof}

\subsection{Valuations at powers of $p$}

Although the valuation formula in Theorem \ref{main1} does not hold for arbitrary Lubin--Tate formal groups in general, it remains valid when the index is a power of $p$.
We prove this by comparing an arbitrary Lubin--Tate formal group $\G$ attached to $p$ with $\LT$.

Since there exists $\psi(Z) \in Z+Z^2\mathcal{O}_L[[Z]]$ such that $\log_\G=\log_{\LT} \circ \psi$, each $P_n^{\G}$ can be expressed as
\begin{align}\label{sum_of_P_n}
P_n^{\G}(X)=P_n^{\LT}(X)+\sum_{i=1}^{n-1} c_i P_i^{\LT}(X)
\end{align}
for some $c_1, \cdots, c_{n-1} \in \mathcal{O}_L$.
Furthermore, if we fix a generator $G_{\LT}$ of $\operatorname{Hom}_{\mathcal{O}_{\C_p}}(\LT, \widehat{\mathbb{G}_m})$, then
\begin{align*}
(G_{\LT} \circ \psi)(Z)=\exp(\Omega_{\LT} \log_\G(Z))-1
\end{align*}
is a generator of $\operatorname{Hom}_{\mathcal{O}_{\C_p}}(\G, \widehat{\mathbb{G}_m})$ and hence we may assume that $\Omega_\G=\Omega_{\LT}$.
From these observations, we obtain the following proposition.

\begin{prop}\label{thm_for_powers}
For any $n \geq 0$, we have
\begin{align*}
v_p(P_{p^n}^\G(\Omega_\G))=w_f(p^n).
\end{align*}
%
\end{prop}

\begin{proof}
We first note that the claim for $\LT$ follows from Theorem \ref{main1}.

For a general $\G$, \eqref{sum_of_P_n} gives
\begin{align*}
P_{p^n}^{\G}(\Omega_{\G})=P_{p^n}^{\LT}(\Omega_{\LT})+\sum_{i=1}^{p^n-1} c_i P_i^{\LT}(\Omega_{\LT})
\end{align*}
for some $c_1, \cdots, c_{p^n-1} \in \mathcal{O}_L$, since we take $\Omega_\G=\Omega_{\LT}$.
Since Theorem \ref{main1} and Lemma \ref{diff_p} imply that
\begin{align*}
v_p(P_i^{\LT}(\Omega_{\LT}))=w_f(i)>w_f(p^n)
\end{align*}
for any $1 \leq i \leq p^n-1$, we obtain
\begin{align*}
v_p(P_{p^n}^{\G}(\Omega_{\G}))=v_p(P_{p^n}^{\LT}(\Omega_{\LT}))=w_f(p^n).
\end{align*}
\end{proof}

\section{Applications to $p$-adic functional analysis}

In this section, we prove Theorem \ref{cor1} by adapting the arguments of \cite{Ka26}. 
Since most of the arguments carry over without modification, we only give the details of the auxiliary estimates that require the properties of the weighted Monna map $w_f$, and refer to \cite{Ka26} for the remaining arguments.

In the following, we assume that $\G=\LT$ and let $N$ be a positive integer.

\subsection{Auxiliary estimates}

We first recall some basic properties of $w_f$ and $\eta_N$ that will be used in the proof of Theorem \ref{cor1}.
The following two lemmas can be proved by the same arguments as \cite[Lemmas 2.3 and 2.4]{Ka26}, and hence we omit their proofs.

\begin{lem}[cf. {\cite[Lemma 2.3]{Ka26}}]\label{ineq_w_eta}
\begin{enumerate}
\item For any $n, k \geq 0$, we have
\begin{align*}
w_f(p^kn) \leq u \sum_{i=0}^{\infty} n_i.
\end{align*}
\item For any $n \geq 0$, we have $w_f(n) \geq \eta_N(n)$. Moreover, equality holds if and only if $0 \leq n \leq q^{N+1}-1$.
\end{enumerate}
\end{lem}

\begin{lem}[cf. {\cite[Lemma 2.4]{Ka26}}]\label{subadd_eta}
For any $n_1, n_2 \geq 0$, we have 
\begin{align*}
\eta_N(n_1+n_2) \leq \eta_N(n_1)+\eta_N(n_2).
\end{align*}
\end{lem}

We next estimate the coefficients of the local expansions of the Schneider–Teitelbaum polynomials.
As in \cite[Subsection 3.1]{Ka26}, the generating series of the Schneider–Teitelbaum polynomials gives the expansion
\begin{align*}
P_n((a+p^NY)\Omega)=\sum_{d=0}^n \left( \sum_{\substack{\nu \\ |\nu|=d, M(\nu)\leq n}} P_{n-M(\nu)}(a\Omega) \frac{p^{\sum_{i \geq 0}(N-i)\nu_i} \Omega^{d}}{\prod_{i \geq 0}\nu_i!} \right) Y^{d}
\end{align*}
for $a \in \mathcal{O}_L$, where 
\begin{align*}
M(\nu) \coloneqq \sum_{i=0}^{\infty} \nu_iq^i, \qquad |\nu| \coloneqq \sum_{i=0}^\infty \nu_i
\end{align*}
for a sequence $\nu=(\nu_i)_{i \geq 0}$ of nonnegative integers with finite support.
To estimate the $p$-adic valuations of the coefficients, we define
\begin{align*}
\phi_{N, i}(t) \coloneqq (u+N-i)t-v_p(t!)
\end{align*}
for nonnegative integers $i$ and $t$.
Then, for any $\nu=(\nu_i)_{i \geq 0}$ with $|\nu|=d$ and $M(\nu) \leq n$, we have
\begin{align}\label{val of target}
\begin{split}
&v_p \left( P_{n-M(\nu)}(a\Omega) \frac{p^{\sum_{i \geq 0}(N-i)\nu_i} \Omega^{d}}{\prod_{i \geq 0}\nu_i!} \right) \\
=& v_p(P_{n-M(\nu)}(a\Omega)) + \sum_{i \geq 0} v_p \left( \frac{p^{(N-i)\nu_i} \Omega^{\nu_i}}{\nu_i!} \right) \\
\geq & w_f(n-M(\nu)) + \sum_{i \geq 0} \phi_{N, i}(\nu_i)
\end{split}
\end{align}
by Corollary \ref{comm_main1}.

We now establish the auxiliary estimates needed to adapt the proof of \cite{Ka26} to the present setting.

\begin{lem}[cf. {\cite[Lemma 3.2]{Ka26}}]\label{phi_eta}
For any nonnegative integers $t$ and $i$, we have $\phi_{N, i}(t) \geq \eta_N(tq^i)$.
Moreover, equality holds if and only if at least one of the following conditions holds:
\begin{enumerate}
\item $t=0$,
\item $i=N$ and $q \mid t$,
\item $i=N+1$.
\end{enumerate}
\end{lem}

\begin{proof}
We consider the three cases: $0 \leq i \leq N-1$, $i=N$ and $i \geq N+1$.
However, the arguments for the first and third cases are exactly the same as those in the proof of \cite[Lemma 3.2]{Ka26}, and hence we omit them.

We suppose that $i=N$ and write $t=aq+b$ with $a \geq 0$ and $0 \leq b \leq q-1$.
We note that
\begin{align}\label{val_of_t!}
v_p(t!)=(s+1)a+v_p(a!)+v_p(b!)
\end{align}
by direct computation, and that
\begin{align*}
ub-v_p(b!)=w_f(b)
\end{align*}
by Proposition \ref{explicit_w_f}.
Since $tq^N=aq^{N+1}+bq^N$ and $0 \leq bq^N <q^{N+1}$, it follows from the definition of $\eta_N$ and \eqref{val_of_t!} that
\begin{align*}
&\phi_{N, N}(t)-\eta_N(tq^N) \\
=&u(aq+b)-v_p(t!)-\{ w_f(bq^N)-a(1-u)-v_p(a!) \} \\
=&ub-v_p(b!)-w_f(bq^N)+a \{ uq-(s+1)+(1-u) \}.
\end{align*}
Since $uq-(s+1)+(1-u)=u(q-1)-s=0$, we obtain
\begin{align*}
\phi_{N, N}(t)-\eta_N(tq^N)=w_f(b)-w_f(bq^N)=\sum_{i=0}^{N-1} D_f(bq^i) \geq 0
\end{align*}
by Proposition \ref{prop2.1} (2).
Since $D_f(b)>0$ whenever $b>0$, $\phi_{N, N}(t)=\eta_N(tq^N)$ if and only if $b=0$, hence $q$ divides $t$.
This completes the proof of the lemma.
\end{proof}

\subsection{Proofs of Theorem \ref{cor1} and Corollary \ref{cor2}}


We are now ready to prove Theorem \ref{cor1} using the auxiliary estimates established in the previous subsection.

\begin{proof}[Proof of Theorem \ref{cor1}]
By Lemmas \ref{ineq_w_eta}, \ref{subadd_eta} and \ref{phi_eta}, the same argument as in the proof of \cite[Corollary 3.3]{Ka26} gives
\begin{align*}
||P_n(x\Omega)||_N \leq p^{-\eta_N(n)}.
\end{align*}
To prove the reverse inequality, we consider the coefficient of $Y^l$ in the expansion of $P_n((1+p^NY)\Omega)$. 
As in \cite[Subsection 3.2]{Ka26}, the term corresponding to the sequence $\nu'=(\nu'_i)_{i\ge0}$, defined by
\begin{align*}
\nu'_i=
\begin{cases}
l & \text{if }i=N+1, \\
0 & \text{otherwise},
\end{cases}
\end{align*}
has valuation $\eta_N(n)$, while all other terms have strictly larger valuations. 
Indeed, the proof of \cite[Lemma 3.5]{Ka26} carries over without essential changes, using Theorem \ref{main1}, Lemmas \ref{ineq_w_eta}, \ref{subadd_eta} and \ref{phi_eta} in place of the corresponding results in \cite{Ka26}. 
Consequently, we obtain
\begin{align*}
||P_n(x\Omega)||_N \geq p^{-\eta_N(n)}.
\end{align*}
Combining the two inequalities proves the theorem.
\end{proof}

We now prove Corollary \ref{cor2}. 
Recall that Bannai and Kobayashi \cite{BK16} introduced a normalization factor $\underline{\gamma}(k)\in\mathbb{C}_p$ for each nonnegative integer $k$, satisfying
\begin{align*}
|\underline{\gamma}(k)|=\min_{0\leq i\leq k}
\left\{\left|\frac{i!}{\Omega^i}\right|\right\}.
\end{align*}
In particular, the functions
\begin{align*}
e_{N,n}(x):=
\underline{\gamma}\left(\left\lfloor\frac{n}{q^N}\right\rfloor\right)
P_n(x\Omega)
\end{align*}
form a Banach basis of $\mathrm{LA}_N(\mathcal{O}_L,\mathbb C_p)$.

\begin{proof}[Proof of Corollary \ref{cor2}]
The proof is identical to that of \cite[Corollary 1.4]{Ka26}. 
Indeed, the normalization factor satisfies
\begin{align*}
v_p\left( \underline{\gamma}\left(\left\lfloor\frac{n}{q^N}\right\rfloor\right) \right)=v_p\left(\frac{(lq)!}{\Omega^{lq}}\right)=l(1-u)+v_p(l!)
\end{align*}
by \cite[Proposition 3.4]{BK16}.
Hence Theorem \ref{cor1} gives
\begin{align*}
||e_{N,n}(x)||_N=p^{-w_f(r)}.
\end{align*}
Since $v_p(P_r(\Omega))=w_f(r)$ by Theorem \ref{main1}, the assertion follows.
\end{proof}

\vspace{10pt}

\noindent
Institute of Mathematics for Industry, Kyushu University,\\
744, Motooka, Nishi-ku, Fukuoka, 819-0395, Japan,\\
E-mail address: \textbf{yu.katagiri.s3@gmail.com}\\

\end{document}